\documentclass[12pt, reqno,dvipsnames]{amsart}
\usepackage{amssymb,latexsym, mathtools,tikz}
\usepackage[pdfencoding=auto, psdextra]{hyperref}
\usepackage[margin=1in,headheight=13.6pt]{geometry}
\usetikzlibrary{arrows,decorations.markings, cd}
\tikzstyle arrowstyle=[scale=1]
\tikzstyle directed=[postaction={decorate,
decoration={markings,mark=at position .65 with {\arrow[arrowstyle]{stealth}}}}]
\usepackage{subcaption, mathdots, adjustbox}

\usepackage{graphicx}

\usepackage{xcolor}
\usepackage{microtype}
\usepackage{enumerate}
\usepackage{graphicx}
\usepackage{float}
\usepackage{placeins}
\usepackage{mdframed}
\usepackage{amssymb}
\usepackage{esint}
\usepackage{cool}
\usepackage[all,cmtip]{xy}
\usepackage{mathtools}
\usepackage{amstext} 
\usepackage{array}   
\usepackage[shortlabels]{enumitem}
\usepackage{ytableau}
\usepackage[ruled,vlined]{algorithm2e}
\usepackage{mathbbol}

\usepackage{cleveref}
\usepackage[utf8]{inputenc}
\usepackage{comment}

\newcolumntype{L}{>{$}l<{$}} 
\newtheorem{lemma}{Lemma}[section]
\newtheorem{theorem}[lemma]{Theorem}
\newtheorem{prop}[lemma]{Proposition}
\newtheorem{cor}[lemma]{Corollary}

\theoremstyle{remark}
\newtheorem{remark}[lemma]{Remark}

\theoremstyle{definition}
\newtheorem{definition}[lemma]{Definition}
\newtheorem{example}[lemma]{Example}
\newtheorem{setup}[lemma]{Setup}
\newtheorem{notation}[lemma]{Notation} 
\newtheorem{construction}[lemma]{Construction} 
\newtheorem{obs}[lemma]{Observation}
\newtheorem{chunk}[lemma]{}

\newcommand{\one}{\mathbb{1}}
\newcommand{\zero}{\mathbb{0}}

\newcommand{\inn}{\operatorname{in}}
\newcommand{\w}{\wedge}

\newcommand{\LS}{\mathrm{LS}}
\newcommand{\tabs}{\mathrm{tab}}

\newcommand{\coker}{\operatorname{coker}}

\newcommand{\tor}{\operatorname{Tor}}

\newcommand{\im}{\operatorname{im}}

\newcommand{\sk}{\operatorname{sk}}

\newcommand{\Supp}{\operatorname{Supp}}

\newcommand{\conv}{\operatorname{conv}}

\newcommand{\rank}{\operatorname{rank}}

\newcommand{\supp}{\operatorname{supp}}

\newcommand{\mdeg}{\operatorname{mdeg}}

\newcommand{\kos}{\textrm{kos}}
\newcommand{\first}{\operatorname{first}}

\newcommand{\ff}{\mathbf f}
\newcommand{\hh}{\mathbf h}

\renewcommand{\aa}{\mathbf a}
\newcommand{\bb}{\mathbf b}
\newcommand{\cc}{\mathbf c}
\newcommand{\dd}{\mathbf d}

\newcommand{\mm}{\mathbf m}
\newcommand{\uu}{\mathbf u}

\newcommand{\ww}{\mathbf w}
\newcommand{\xx}{\mathbf x}

\newcommand{\cB}{\mathcal{B}}
\newcommand{\cC}{\mathcal{C}}

\newcommand{\cG}{\mathcal{G}}
\newcommand{\cH}{\mathcal{H}}

\newcommand{\cP}{\mathcal{P}}

\newcommand{\FF}{\mathbb{F}}
\newcommand{\GG}{\mathbb{G}}

\newcommand{\NN}{\mathbb{N}}

\newcommand{\RR}{\mathbb{R}}

\newcommand{\UU}{\mathbb{U}}		

\newcommand{\ZZ}{\mathbb{Z}}

\newcommand{\m}{\mathfrak{m}}

\newcommand{\Xv}{\check{X}}

\newcommand{\ra}{\rightarrow}

\DeclarePairedDelimiter\abs{\lvert}{\rvert}%
\DeclareMathOperator{\lcm}{lcm}

\begin{document}

\title[Polarizations and Hook Partitions]{Polarizations and Hook Partitions}

\author{Ayah Almousa}
\address{University of Minnesota - Twin Cities}
\email{almou007@umn.edu}
\urladdr{\url{http://umn.edu/~almou007}}

\author{Keller VandeBogert}
\address{University of Notre Dame}
\email{kvandebo@nd.edu}

\urladdr{\url{https://sites.google.com/view/kellervandebogert/}}

\date{\today}
\begin{abstract}
In this paper, we relate combinatorial conditions for polarizations of powers of the graded maximal ideal with rank conditions on submodules generated by collections of Young tableaux. We apply discrete Morse theory to the hypersimplex resolution introduced by Batzies--Welker to show that the $L$-complex of Buchsbaum and Eisenbud for powers of the graded maximal ideal is supported on a CW-complex. We then translate the ``spanning tree condition'' of Almousa--Fl\o ystad--Lohne characterizing polarizations of powers of the graded maximal ideal into a condition about which sets of hook tableaux span a certain Schur module. As an application, we give a complete combinatorial characterization of polarizations of so-called ``restricted powers'' of the graded maximal ideal.
\end{abstract}
\maketitle

\section{Introduction}

Let $S = k[x_1 , \dots , x_n]$ be a standard graded polynomial ring over a field $k$. The Taylor resolution (introduced in \cite{taylor66}) is well-known to be a free resolution of any ideal generated by monomials in $S$; this complex has many convenient properties, not least of which is the fact that it is locally isomorphic to an exterior algebra, implying that the Taylor resolution behaves in a manner that is very similar to the Koszul complex. In \cite{bayer1998monomial}, Bayer, Peeva, and Sturmfels observed that this implies that every monomial ideal has a free resolution supported on a simplicial complex. It is clear that not every monomial ideal has a \emph{minimal} free resolution supported on a simplicial complex, and indeed, work of Velasco \cite{velasco2008minimal} has shown that even the much weaker notion of CW-complexes is not general enough to support all resolutions of monomial ideals. 

One natural question that arises from the above considerations is the following: given that the minimal free resolution of every monomial ideal $I$ is a direct summand of the Taylor complex $T_\bullet$, how can one extract a strictly smaller subcomplex $F_\bullet \subset T_\bullet$ satisfying
\begin{enumerate}
    \item $F_\bullet$ is still supported on a cell complex, and
    \item $F_\bullet$ is a free resolution of $I$?
\end{enumerate}
Batzies and Welker give one possible answer to this question using discrete Morse theory; the basic idea is as follows: suppose that $X$ is a cell complex supporting the resolution of some ideal $I$. If the associated graph (see Construction \ref{const: DMT}) admits an acyclic matching on some edge set $A$, then one can construct an associated Morse complex $X_A$ that remains acyclic, is closer to being minimal, and is also supported on a cell complex. Let $\m = (x_1 , \dots , x_n) \subset S$. One application of these techniques by Batzies is a proof that the Eliahou-Kervaire resolution for resolving $\m^d$ is supported on a cell complex (it was later proved that the Eliahou-Kervaire in full generality is cellular independently by Mermin \cite{mermin2010eliaiiou} and Clark \cite{clark2012minimal}).

As it turns out, there is another well-known complex providing a minimal free resolution of powers of the maximal ideal, introduced by Buchsbaum and Eisenbud in \cite{buchsbaum1975generic}. The free modules building these $L$-complexes are Schur modules corresponding to appropriate hook partitions, whose basis elements are represented by semistandard tableaux subject to so-called straightening relations. We show that the $L$-complexes of Buchsbaum and Eisenbud are also supported on a cell complex using techniques similar to those of Batzies; one step in this proof is the observation that there is a simple bijection between the basis elements of the hypersimplex resolution (see Definition \ref{def: hypersimplex}) and basis elements of an associated enveloping algebra. 

Batzies' hypersimplex resolution in particular keeps track of all possible linear syzygies on the monomial minimal generating set of $\m^d$, which can be encoded as a graph. In \cite{polarizations}, Almousa, Fl\o ystad, and Lohne show that this graph of linear syzygies can be used to characterize all possible \emph{polarizations} (see Definition \ref{def: polarization}) of $\m^d$. We use the aforementioned bijection of basis elements to translate conditions on spanning trees contained within the graph of linear syzygies to rank conditions on submodules generated by elements of an appropriate Schur module. This yields a dictionary between the notation and terminology introduced in \cite{polarizations} and well-established notions arising in the context of Schur modules. Moreover, we extend the results of \cite{polarizations} to polarizations of \emph{restricted powers} of the maximal ideal, and give an explicit algorithm for checking whether a given graph of linear syzygies induces a well-defined isotone map. This algorithm opens the door to methods of computing all possible polarizations of $\m^d$ for any number of variables.

The paper is organized as follows. In Section \ref{sec: dmt}, we recall the definition of a cellular resolution and summarize the results of Batzies--Welker \cite{batzies2002discrete} on discrete Morse theory for cellular resolutions. In Section \ref{sec: Lcomplex}, we recall the construction for the $L$-complex of Buchsbaum and Eisenbud, which is a minimal free resolution of powers of complete intersections. In Section \ref{sec: hypersimplex}, we begin our study of the so-called \emph{hypersimplicial complex} $\cH_n^d$ (Definition \ref{def: hypersimplex}) of Batzies--Welker. We observe that the cells $\cH_n^d$ correspond to hook tableaux, and that $\cH_n^d$ supports a free resolution of a power of the graded maximal ideal coming from a certain double complex (see \ref{chunk: hypersimplexHooks} and Figure \ref{fig: doubleCplx}). With this perspective in mind, we apply discrete Morse theory in a novel way to obtain a CW-complex supporting the $L$-complex of Buchsbaum and Eisenbud for powers of the graded maximal ideal (Propositions \ref{prop: DMTres} and \ref{prop: LcplxCellular}). 

In Section \ref{sec: pols}, we summarize the results from \cite{polarizations} giving a complete characterization of all polarizations of powers of the graded maximal ideal. In Section \ref{sec: dictionary}, we translate the machinery of Section \ref{sec: pols} to the language of hook tableaux (Propositon \ref{prop: dictionary}). We utilize this perspective to give a novel characterization of polarizations of powers of the graded maximal ideal $\mm^d$ in terms of hook tableaux which span the Schur module $L^2_d$ (see Theorem \ref{thm: spanningTree}). 

In Section \ref{sec: boundedGens}, we extend the results from Sections \ref{sec: hypersimplex} and \ref{sec: dictionary} from powers of the graded maximal ideal to a larger class of ideals called \emph{restricted powers} of the graded maximal ideal. In particular, we give a complete characterization of all polarizations of such ideals in Theorem \ref{thm: generalizedPols}, which is a direct extension of the results on polarizations of $\m^d$ in \cite{polarizations}.

\section{Frames and Discrete Morse Theory for Cellular Resolutions}\label{sec: dmt}

In this section, we recall some important notions on cellular resolutions and Discrete Morse theory for cellular resolutions. For further exposition on frames and cellular resolutions, we refer the reader to \cite{peeva2010graded} and \cite{peeva2011frames}. We will use the terminology of frames to give a convenient framework (pun intended) for defining cellular resolutions. Proposition \ref{prop: DMTres} will be essential for proving that the $L$-complexes of Buchsbaum and Eisenbud are cellular. We begin by adopting the following setup:

\begin{setup}\label{set: frames} Let $S = k[x_1,\dots, x_n]$ be a polynomial ring over a field $k$. Let $M$ be a monomial ideal in $S$ minimally generated by monomial $m_1,\dots, m_r$. Let $L_M$ denote the set of least common multiples of subsets of $m_1,\dots, m_r$. By convention, $1\in L_M$ is considered to be the lcm of the empty set.
\end{setup}

\begin{definition}\label{def: frame} Adopt notation and hypotheses of Setup \ref{set: frames}. A \emph{frame} (or an $r$-frame) $\UU$ is a complex of finite $k$-vector spaces with differential $\partial$ and a fixed basis that satisfies the following conditions:
\begin{enumerate}[(1)]
    \item $U_i = 0$ for $i<0$ and $i\gg 0$,
    \item $U_0 = k$
    \item $U_1 = k^r$
    \item $\partial(w_j) = 1$ for each basis vector $w_j$ in $U_1 = k^r$.
\end{enumerate}
\end{definition}

Given a complex of free modules over some polynomial ring, it is easy to obtain a frame by setting all variables equal to $1$. Conversely, given a frame $\UU$, one may construct a multigraded complex $\GG$ of finitely generated free multigraded $S$-modules with multidegrees in $L_M$ using the following construction due to Peeva and Velasco \cite{peeva2011frames}. 

\begin{construction}\label{const: homogenization} Adopt notation and hypotheses of Setup \ref{set: frames}. Let $\UU$ be an $r$-frame. Set
$$
G_0 = S \text{ and } G_1 = S(-m_1)\oplus \dots \oplus S(-m_r).
$$
Let $\overline{v}_1,\dots, \overline{v}_p$ and $\overline{u}_1,\dots, \overline{u}_q$ be the given bases of $U_i$ and $U_{i-1}$, respectively. Let $u_1,\dots, u_q$ be the basis of $G_{i-1} = S^q$ chosen on the previous step of the induction. Introduce $v_1,\dots, v_p$ that will be a basis of $G_i = S^p$. If
$$
\partial(\overline{v}_j) = \sum_{1\leq s\leq q} \alpha_{sj} \overline{u}_s
$$
with coefficients $\alpha_{sj} \in k$, then set
\begin{gather*}
    \mdeg(v_j) = \lcm(\mdeg(u_s)\mid \alpha_{sj}\neq 0)\\
    G_i = \bigoplus_{1\leq j\leq p} S(-\mdeg(v_j))\\
    d(v_j) = \sum_{1\leq s \leq q} \alpha_{sj} \frac{\mdeg(v_j)}{\mdeg(u_s)}\cdot u_s.
\end{gather*}
Note that given a monomial $m = \xx^\bb$ where $\bb\in \NN^n$, in our notation $\mdeg(m)$ is equal to the monomial $m$ itself, rather than the $\NN^n$ degree $\bb$. 
Clearly $\coker(d_1) = S/M$ and the differential $d$ is homogeneous by construction. Call $\GG$ the \emph{$M$-homogenization} of $\UU$.
\end{construction}

The following simple criterion by Peeva and Velasco \cite{peeva2011frames} determines when a frame supports a graded free resolution of $S/M$. The abridged version of this result states that exactness can be checked by only considering multihomogeneous strands.

\begin{prop}\label{prop: homogRes} The sequence of modules and homomorphisms $\GG$ as in Construction \ref{const: homogenization} is a complex. Moreover, if $\GG(\leq m)$ is the subcomplex of $\GG$ generated by the multihomogeneous basis elements of multidegrees dividing $m$, then $\GG$ is a free multigraded resolution of $S/M$ if and only if for all monomials $1\neq m\in L_M$, the frame of the complex $G(\leq m)$ is exact.
\end{prop}

A natural source of frames that can be used to support resolutions of monomial ideals are provided by CW-complexes, since the conditions $(1)-(4)$ of Definition \ref{def: frame} are trivially satisfied. 

\begin{notation}\label{not: cellRes}  Let $X$ be a regular CW-complex, and denote  by $X^{(i)}$ the set of $i$-cells of $X$ and by $X^{(\ast)} \coloneqq \bigcup_{i\geq 0} X^{(i)}$ the set of all cells of $X$.
Denote by $\widetilde \cC(X;k)$ the augmented oriented reduced cellular chain complex of $X$ over $k$ with
$$
\widetilde \cC(X;k)_i = \bigoplus_{c\in X^{(i)}} k e_c
$$
where $e_c$ denotes the basis element corresponding to the face $c\in X^{(i)}$, and the differential $\partial$ acts as
$$
\partial(e_c) = \sum_{c\geq c'\in X^{(i-1)}} [c,c'] e_{c'}
$$
where $[c,c']$ is the coefficient in the differential of the cellular homology of $X$.
\end{notation}

With the above notation in mind, we use the language of frames to define cellular resolutions.

\begin{definition}\label{def: cellRes}
Adopt notation of Notation \ref{not: cellRes}. Assume that $\abs{X^{(0)}} = r$ and $M = (m_1, \dots, m_r)$ is a monomial ideal in a polynomial ring $S$. Label each $0$-cell of $X$ by a minimal generator $m_i$ of $M$. After shifting $\widetilde \cC(X;k)$ in homological degree, $\widetilde \cC(X;k)[-1]$ is a frame. Denote by $\FF_X$ the $M$-homogenization of $\widetilde \cC(X;k)$ as in Construction \ref{const: homogenization}. The complex $\FF_X$ is \emph{supported on} $X$. The complex $\FF_X$ is a \emph{cellular resolution} if it is exact.
\end{definition}

\begin{definition}\label{def: CWrestrict} Let $M = (m_1,\dots, m_r)$ be a monomial ideal in a polynomial ring $S$, and let $X$ be a regular CW-complex with $0$-cells labeled by the generators of $M$. 
The multidegree of each vertex of $X$ is given by its monomial label. Define a face $c$ to have multidegree
$$
\mdeg(c) = \lcm(m_i \mid m_i\in c).
$$
By convention, $\mdeg(\emptyset) = 1$.
Define the following subcomplexes of $X$:
\begin{align*}
    X_{\leq m} &\coloneqq \{c\in X \mid \mdeg(c) \text{ divides } m\}\\
    X_{< m} &\coloneqq \{c\in X \mid \mdeg(c) \text{ strictly divides } m\}.
\end{align*}
\end{definition}

The following Proposition is an immediate consequence of Proposition \ref{prop: homogRes} combined with the notation and hypotheses introduced in Definition \ref{def: CWrestrict}.

\begin{prop}\label{prop: restrictionAcyclic} Let $M = (m_1,\dots, m_r)$ be a monomial ideal in a polynomial ring $S$, and let $X$ be a regular $CW$-complex with $0$-cells labeled by the minimal generators of $M$. The complex $\FF_X$ from Definition \ref{def: cellRes} is a free resolution of $S/M$ if and only if for all multidegrees $1\neq m \in L_M$, the complex $X_{\leq m}$ is acyclic over $k$.
\end{prop}

Next, we introduce some of the basic machinery of discrete Morse theory for cellular resolutions. Discrete Morse theory was developed by Forman in \cite{forman1998morse} to extend the ideas from Morse theory in differential geometry to CW complexes. The interested reader is encouraged to consult Forman's survey paper \cite{forman2002user} for further reading on discrete Morse theory. 
The application of discrete Morse theory to the study of cellular resolutions was first explored by Batzies and Welker in \cite{batzies2002discrete} as a method of ``cutting down" a large cellular resolution in such a manner that the resulting subcomplex is also a cellular resolution.

\begin{construction}\label{const: DMT} Adopt Notation \ref{not: cellRes}. Let $G_X$ be the directed graph on the set of cells of $X$ whose set $E_X$ of edges is given by $c\ra c'$ for $c'\subseteq c$ and $\dim(c') = \dim(c)-1$. A \emph{discrete Morse function} arises from a set $A\subseteq E_X$ of edges in $G_X$ satisfying:
\begin{enumerate}
    \item each cell occurs in at most one edge of $A$, and
    \item the graph $G_X^A$ with edge set
$$
E_X^A \coloneqq (E_X\setminus A) \cup \{c'\ra c\mid c\ra c' \in A\}
$$
is acyclic (i.e., it does not contain a directed cycle).
\end{enumerate} 
Such a set $A\subseteq E_X$ is called an \emph{acyclic matrching} of $G_X$. A cell of $X$ is \emph{$A$-critical} with respect to $A$ if it is not contained in any edge of $A$. An acyclic matching is \emph{homogeneous} if $c\ra c'\in A$ implies that $\mdeg(c) = \mdeg(c')$.
\end{construction}

The proof of the following proposition can be found in the appendix of \cite{batzies2002discrete}, and shows that acyclic matchings can be used to induce acyclic subcomplexes that are also supported on cell complexes.

\begin{prop}\label{prop: DMTres} Let $X$ be a regular CW-complex which supports a free resolution of a monomial ideal $M$, and let $A$ be a homogeneous acyclic matching of $G_X$. Then there is a (not necessarily regular) CW-complex $X_A$ whose $i$-cells are in one-to-one correspondence with the $A$-critical $i$-cells of $X$ such that $X_A$ is homotopy equivalent to $X$. 

Moreover, $X_A$ inherits a multigrading from $X$, and for any multidegree $\alpha$ and restriction $A_{\leq \alpha}$ of $A$ to $X_{\leq \alpha}$, one has
$$
X_{\leq \alpha} \simeq (X_{\leq \alpha})_{A_{\leq \alpha}} \cong (X_A)_{\leq \alpha}.
$$
In particular, $X_A$ also supports a cellular resolution of the ideal $M$.
\end{prop}

\begin{definition}\label{def: morseCx}
The complex $X_A$ of Proposition \ref{prop: DMTres} is called the \emph{Morse complex} of $X$ for the matching $A$. 
\end{definition}

\begin{remark}
The explicit construction of the Morse complex $X_A$ from an acyclic matching $A$ is quite technical, and will not be included in the current paper. The interested reader is encouraged to consult the appendix of \cite{batzies2002discrete} for more details.
\end{remark}

\section{Background on $L$-complexes}\label{sec: Lcomplex}

The material up until Proposition \ref{prop:resnofpower}, along with proofs, can be found in \cite{buchsbaum1975generic} or Section $2$ of \cite{el2014artinian}. The goal of this section is to a give a brief jog through the $L$-complexes of Buchsbaum and Eisenbud and to make clear our conventions on Young tableaux. For further details on Schur modules and their use in the construction of free resolutions, one may consult Weyman's book \cite{weyman2003}. The following notation will be in play for the remainder of this section.

\begin{notation}
Let $R$ be a polynomial ring over a field $k$. Let $F$ be a free $R$-module of rank $n$ with basis $f_1,\dots, f_n$. Denote by $S_d(F)$ the $d$th symmetric power of $F$, and by $\bigwedge^n F$ the $n$th exterior power of $F$. Let $J = \{ j_1 < j_2 < \dots < j_k\} \subset [n]$. Define
$$
f_J \coloneqq f_{j_1} \wedge \dots \wedge f_{j_k}\in \bigwedge^k F.
$$
If $\alpha = (\alpha_1, \dots, \alpha_n) \in \NN^n$ such that $\sum_i \alpha_i = d$, set
$$
f^\alpha \coloneqq f_1^{\alpha_1} f_2^{\alpha_2} \dots f_n^{\alpha_n} \in S_d(F).
$$
\end{notation}

\begin{setup}\label{set:Lcomplexsetup}
Let $F$ denote a free $R$-module of rank $n$, and $S = S(F)$ the symmetric algebra on $F$ with the standard grading. Define a complex
$$\xymatrix{\cdots \ar[r] & \bigwedge^{a+1} F \otimes_R S_{b-1} \ar[r]^-{\kappa_{a+1,b-1}} & \bigwedge^{a} F \otimes_R S_{b} \ar[r]^-{\kappa_{a,b}} & \cdots}$$
where the maps $\kappa_{a,b}$ are defined as the composition
\begin{equation*}
    \begin{split}
         \bigwedge^{a} F \otimes_R S_{b} &\to \bigwedge^{a-1} F \otimes_R F \otimes_R S_{b} \\
         & \to \bigwedge^{a-1} F \otimes_R S_{b+1}
    \end{split}
\end{equation*}
where the first map is comultiplication in the exterior algebra and the second map is the standard module action (where we identify $F = S_1 (F)$).
Define
$$L_b^a (F) := \ker \kappa_{a,b}.$$
Let $\psi: F \to R$ be a morphism of $R$-modules with $\im (\psi)$ an ideal of grade $n$. Let $\kos^\psi : \bigwedge^i F \to \bigwedge^{i-1} F$ denote the standard Koszul differential; that is, the composition
\begin{equation*}
    \begin{split}
        \bigwedge^i F &\to F \otimes_R \bigwedge^{i-1} F  \quad \textrm{(comultiplication)} \\
        &\xrightarrow[]{\psi\otimes 1} R\otimes_R \bigwedge^{i-1}F \cong \bigwedge^{i-1} F \quad \textrm{(module action).} \\
    \end{split}
\end{equation*}
Explicitly, if $J = \{j_1 < \dots < j_k\}$, then
$$
\kos^\psi(f_J) = \sum_{i\in [k]} (-1)^i \psi(f_{j_i})\cdot f_{J\setminus j_i}.
$$
\end{setup}

\begin{definition}\label{def:Lcomplexes}
Adopt notation and hypotheses of Setup \ref{set:Lcomplexsetup}. Define the complex \textit{$L$-complex} to be the complex
\begin{equation*}
    \begin{split}
        &L(\psi , b) : \xymatrix{0 \ar[r] & L_b^{n-1} \ar[rr]^-{\kos^\psi \otimes 1} & & \cdots \ar[rr]^{\kos^\psi \otimes 1} & & L_b^0 \ar[r]^-{S_b (\psi)} & R \ar[r] & 0 } \\
    \end{split}
\end{equation*}
where $\kos^\psi \otimes 1  : L_b^a (F) \to L_b^{a-1}$ is induced by making the following diagram commute:
$$\xymatrix{\bigwedge^a F\otimes S_b (F) \ar[rr]^-{\kos^\psi \otimes 1}  & & \bigwedge^{a-1} F \otimes S_b(F)  \\
L_b^a (F) \ar[rr]^-{\kos^\psi \otimes 1} \ar[u] & & L_b^{a-1} (F) \ar[u] \\}$$
\end{definition}

The following Proposition shows that the $L$-complexes constitute a minimal free resolution of powers of complete intersections in general.

\begin{prop}\label{prop:resnofpower}
Let $\psi: F \to R$ be an $R$-module homomorphism from a free module $F$ of rank $n$ such that the image  $\im (\psi)$ is an ideal of grade $n$. Then the complex $L(\psi ,b)$ of Definition \ref{def:Lcomplexes} is a minimal free resolution of $R/\im (\psi)^b$
\end{prop}
We also have (see Proposition $2.5 (c)$ of \cite{buchsbaum1975generic}, or just use Proposition \ref{prop:standardbasis})
\begin{equation*}
    \begin{split}
        &\rank_R L_b^a (F) = \binom{n+b-1}{a+b} \binom{a+b-1}{a}. \\
    \end{split}
\end{equation*}
Moreover, using the notation and language of Chapter $2$ of \cite{weyman2003}, $L_b^a (F)$ is the Schur module $L_{(a+1,1^{b-1})} (F)$. This allows us to identify a standard basis for such modules.

\begin{notation}
We use the English convention for partition diagrams. That is, the partition $(3,2,2)$ corresponds to the diagram
$$ 
\begin{ytableau}
 \ & \ & \ \\
\ & \ \\
\ & \ \\
\end{ytableau}.$$
A Young tableau is standard if it is strictly increasing in both the columns and rows. It is semistandard if it is strictly increasing in the columns and nondecreasing in the rows.
\end{notation}

\begin{prop}\label{prop:standardbasis}
Adopt notation and hypotheses as in Setup \ref{set:Lcomplexsetup}. Then a basis for $L_b^a (F)$ is represented by all Young tableaux of the form
$$\ytableausetup
{boxsize=2em}
\begin{ytableau}
i_0 & j_1 &\cdots & j_{b-1} \\
i_1 \\
\vdots \\
i_a \\
\end{ytableau}$$
with $i_0 < \cdots < i_{a}$ and $i_0 \leq j_1 \leq \cdots \leq j_{b-1}$. 
\end{prop}

\begin{proof}
See Proposition $2.1.4$ of \cite{weyman2003} for a more general statement.
\end{proof}

\begin{remark}\label{rk:rkonstdbasis}
Adopt notation and hypotheses of Setup \ref{set:Lcomplexsetup}. Let $F$ have basis $f_1, \dotsc, f_n$. In the statement of Proposition \ref{prop:standardbasis}, we think of the tableau as representing the element
$$\kappa_{a+1,b-1} (f_{i_0} \w \cdots \w f_{i_{a}} \otimes f_{j_1} \cdots f_{j_{b-1}}) \in \bigwedge^a F \otimes S_b (F).$$
We will often write $f_{i_0} \w \cdots \w f_{i_{a}} \otimes f_{j_1} \cdots f_{j_{b-1}} \in L_b^a (F)$, with the understanding that we are identifying $L_b^a (F)$ with the cokernel of $\kappa_{a+2,b-2} : \bigwedge^{a+2} F \otimes S_{b-2} (F) \to \bigwedge^{a+1} F \otimes S_{b-1} (F)$.
\end{remark}

The following Observation is sometimes referred to as the \emph{shuffling} or \emph{straightening} relations satisfied by tableaux in the Schur module $L_b^a (F)$.

\begin{obs}\label{obs:straighteningrelations}
Any tableau of the form
$$T= \ytableausetup
{boxsize=2em}
\begin{ytableau}
i_0 & j_1 &\cdots & j_{b-1} \\
i_1 \\
\vdots \\
i_a \\
\end{ytableau}$$
with $j_1 \leq \cdots \leq j_{b-1}$ viewed as an element in $L_b^a (F)$ with $b\geq 2$ may be rewritten as a linear combination of other tableaux in the following way:
$$T = \sum_{k=0}^a (-1)^k \ytableausetup
{boxsize=1.7em}
\begin{ytableau}
j_1 & i_k & j_2 &\cdots & j_{b-1} \\
i_0 \\
\vdots \\
\widehat{i_k} \\
\vdots \\
i_a \\
\end{ytableau}.$$
Notice that if $i_0 > j_1$ and $i_0 < \cdots < i_a$, then this rewrites $T$ as a linear combination of semistandard tableaux after reordering the row into ascending order.
\end{obs}

\section{The $L$-complex is cellular}\label{sec: hypersimplex}

In this section, we apply discrete Morse theory to the so-called \emph{hypersimplex resolution} (see Definition \ref{def: hypersimplex}) of $\m^d$ in a novel way to obtain a CW-complex which supports the $L$-complex with the exact basis elements described in Proposition \ref{prop:standardbasis}. In particular, Proposition \ref{prop: LcplxCellular} implies that the $L$-complex of Buchsbaum and Eisenbud is CW-cellular. While Batzies and Welker had previously obtained a minimal cellular resolution of $\m^d$ by finding an acyclic matching on the hypersimplex resolution in \cite{batzies2002discrete}, the minimal resolution they obtain is instead isomorphic to the Eliahou-Kervaire resolution of $\m^d$.

\begin{notation}\label{not: dilatedSimplex}
The notation $\Delta(n,d)$ will denote the dilated $(n-1)$-simplex $d\cdot\Delta^{n-1}$; that is,
$$
\Delta(n,d) = d\cdot \Delta^{n-1} := \left\{(y_1,\dots, y_n)\in \RR^n\mid \sum_{i=1}^n y_i = d, y_i\geq 0 \text{ for } i=1,\dots, n\right\}.
$$
Moreover, the notation $\NN = \{0,1,\dots\}$ will denote the set of nonnegative integers.
\end{notation}

\begin{definition}\label{def: hypersimplex} Let $\cH^d_n$ be the polytopal CW-complex with the underlying space $\Delta(n,d)$, with CW-complex stucture induced by intersection with the cubical CW-complex structure on $\RR^n$ given by the integer lattice $\ZZ^n$. That is, the closed cells of $\cH^d_n$ are given by all hypersimplices
\begin{align*}
C_{\aa, J} &\coloneqq \Delta(n,d) \cap \left\{(y_1,\dots, y_n)\in \RR^n \mid y_i = a_i \text{ for } i\in [n]\setminus J, \text{ and } y_j\in [a_j,a_j+1] \text{ for } j\in J \right\} \\
 &= \conv\left(\aa + \sum_{j\in J} \ell_j \epsilon_j \mid \ell_j \in \{0,1\}, \sum_{j\in J} \ell_j = d - \abs{\aa}\right) 
\end{align*}
with $\aa\in \NN^n$, $J\subset [n]$, $\abs{\aa} \coloneqq \sum_{i\in [n]} a_i$, $\epsilon_i$ the $i$th unit vector in $\RR^n$, either subject to the conditions $\abs{a} = d$ and $J = \emptyset$ (these are the $0$-cells), or the condition $1\leq d -\abs{\aa}\leq \abs{J}-1$.
The CW-complex $\cH^d_n$ is multigraded by setting $\lcm(C_{\aa,J}) \coloneqq \aa + \sum_{j\in J} \epsilon_j$. Call $\cH^d_n$ the \emph{hypersimplicial complex}.

Let $J = (j_0 < \cdots < j_r)$ and use the notation $J_v \coloneqq J\setminus \{j_v\}$. Then the differential of $\cH^d_n$ is given by
\begin{align}\label{eq: hypersimplicialDiffs}
    \partial(C_{\aa,J}) = \begin{cases}
    \sum_{v = 0}^r (-1)^v(C_{\aa,J_v} - C_{\aa+\epsilon_v, J_v}) &\text{ if } 2\leq d-\abs{\aa}\leq \abs{J}-2\\ \\
    \sum_{v = 0}^r (-1)^v C_{\aa, J_v} &\text{ if } 1 = d-\abs{\aa} \leq \abs{J}-2\\ \\
    \sum_{v= 0}^r (-1)^v C_{\aa+\epsilon_j, J_v} &\text{ if } 2\leq d - \abs{\aa} = \abs{J}-1\\ \\
    C_{\aa+e_{j_0},\emptyset} - C_{a+e_{j_1},\emptyset} &\text{ if } d-\abs{\aa} = 1\\ \\
    0 &\text{ if } \abs{\aa} = d.
    \end{cases}
\end{align}
\end{definition}

\begin{chunk}\label{chunk: hypersimplexHooks}
Adopt notation and hypotheses of Setup \ref{set:Lcomplexsetup}. If $J = ( j_0 < j_1 < \dots < j_r )$ and $\aa\in\NN^n$, note that every $r$-dimensional cell $C_{\aa,J}$ corresponds to the element $f_J \otimes f^\aa$ in $\bigwedge^{r+1} F \otimes S_{\abs{a}}(F)$. These elements, in turn, can be represented as hook tableaux with strictly increasing columns and weakly increasing rows:
$$
C_{\aa, J} \quad \longleftrightarrow \quad f_J \otimes f^\aa \quad \longleftrightarrow \quad
\begin{ytableau} j_0 & 1^{a_1} &\dots & n^{a_n} \\
j_1 \\
j_2 \\
\vdots \\
j_r
\end{ytableau}.
$$
We will implicitly use this correspondence to refer to cells $C_{\aa,J}$ of $\cH^d_n$ as tableaux or elements of $\bigwedge^{r+1} F \otimes S_{\abs{a}}(F)$. The $\m^d$-homogenization (see Construction \ref{const: homogenization}) of $\cH^d_n[-1]$ therefore corresponds to the double complex in Figure \ref{fig: doubleCplx} where the maps $\kappa$ and $\kos^\psi\otimes 1$ are as defined in Setup \ref{set:Lcomplexsetup}.
\begin{figure}[ht]
\adjustbox{scale = 0.6}{
\begin{tikzcd}
	&&&&& 0 \\
	&&&& 0 & {\bigwedge^n F} \\
	&&& \iddots & \vdots & \vdots \\
	&& 0 & \dots & {\bigwedge^5 F\otimes S_{d-3}} & {\bigwedge^4 F\otimes S_{d-3}(F)} \\
	& 0 & {\bigwedge^{n} F\otimes S_{d-2}(F)} & \dots & {\bigwedge^4 F \otimes S_{d-2}(F)} & {\bigwedge^3 F\otimes S_{d-2} (F)} \\
	0 & {\bigwedge^n F \otimes S_{d-1}(F)} & \dots & {\bigwedge^4 F\otimes S_{d-1}(F)} & {\bigwedge^3 F\otimes S_{d-1}(F)} & {\bigwedge^2 F\otimes S_{d-1}(F)} \\
	&&&& {} & {F\otimes S_d(F)} & {S_d(F)} & 0
	\arrow["{\kos^\psi\otimes 1}", from=7-6, to=7-7]
	\arrow["\kappa", from=6-6, to=7-6]
	\arrow[bend left = 10, from=6-6, to=7-7]
	\arrow["{\kos^\psi\otimes 1}"', from=6-5, to=6-6]
	\arrow["{\kos^\psi\otimes 1}"', from=6-4, to=6-5]
	\arrow["{\kos^\psi\otimes 1}"', from=6-3, to=6-4]
	\arrow["{\kos^\psi\otimes 1}"', from=6-2, to=6-3]
	\arrow[from=6-1, to=6-2]
	\arrow[from=5-2, to=6-2]
	\arrow[from=5-2, to=5-3]
	\arrow["{\kos^\psi\otimes 1}", from=5-3, to=5-4]
	\arrow["{\kos^\psi\otimes 1}", from=5-4, to=5-5]
	\arrow["{\kos^\psi\otimes 1}", from=5-5, to=5-6]
	\arrow["{\kos^\psi\otimes 1}", from=4-5, to=4-6]
	\arrow["{\kos^\psi\otimes 1}", from=4-4, to=4-5]
	\arrow["\kappa", from=5-6, to=6-6]
	\arrow["\kappa", from=5-5, to=6-5]
	\arrow["\kappa", from=5-3, to=6-3]
	\arrow[from=4-3, to=5-3]
	\arrow["\kappa", from=4-5, to=5-5]
	\arrow["\kappa", from=4-6, to=5-6]
	\arrow["\kappa", from=3-6, to=4-6]
	\arrow["\kappa", from=2-6, to=3-6]
	\arrow[from=1-6, to=2-6]
	\arrow["{\kos^\psi\otimes 1}", from=7-7, to=7-8]
	\arrow[from=4-3, to=4-4]
	\arrow[from=2-5, to=2-6]
	\arrow[from=2-5, to=3-5]
	\arrow["\kappa", from=3-5, to=4-5]
	\arrow["\kappa", from=5-4, to=6-4]
	\arrow["\kappa", from=4-4, to=5-4]
	\arrow["{\kos^\psi\otimes 1}", from=3-5, to=3-6]
	\arrow[from=3-4, to=4-4]
	\arrow[from=3-4, to=3-5]
\end{tikzcd}
}
\caption{Double complex supported on the hypersimplex resolution of a power of the graded maximal ideal.}\label{fig: doubleCplx}
\end{figure}
\end{chunk}

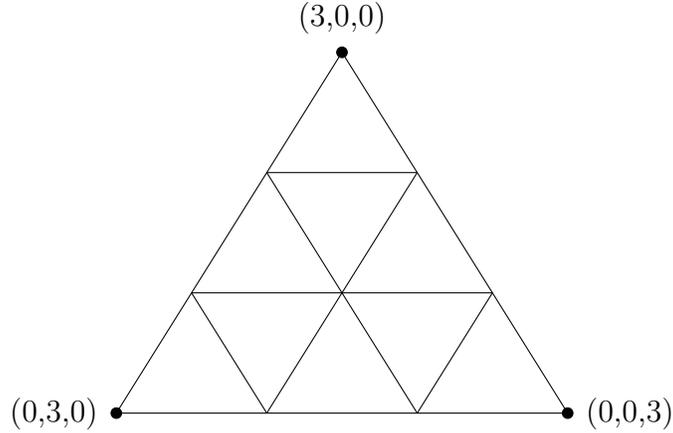
\begin{figure}[ht]
\begin{tikzpicture}
\draw (-3,0)--(3,0);
\draw (-2,1.6)--(2,1.6);
\draw (-1,3.2)--(1,3.2);
\draw (-3,0)--(0,4.8);
\draw (3,0)--(0,4.8) ;
\draw (-1,0)--(1,3.2);
\draw (1,0)--(-1,3.2);
\draw (-1,0)--(-2,1.6);
\draw (1,0)--(2,1.6);
\filldraw[black] (0,4.8) circle (2pt)  node[anchor=south] at (0,4.9){(3,0,0)};
\filldraw[black] (-3,0) circle (2pt)  node[anchor=east] at (-3.1,0){(0,3,0)};
\filldraw[black] (3,0) circle (2pt)  node[anchor=west] at (3.1,0){(0,0,3)};
\end{tikzpicture}
\caption{The hypersimplicial complex $\cH^3_3$.}\label{fig: dilatedSimplex}
\end{figure}

\begin{figure}
    \centering
    \begin{subfigure}[b]{0.47\textwidth}
    \centering
\begin{tikzpicture}[]

\draw    (0, 4.33) -- (-2.5,0)--(2.5,0)--(0,4.33) ;

\draw[very thick, RedOrange] (1.68, 1.44)--(0.85,2.88);

\draw    (-1.68, 1.44)--(0,1.44);

\draw[very thick, BlueGreen] (0,1.44)--(1.68, 1.44) ;

\draw  (-0.85,2.88)--(0,1.44);
\draw[very thick, Plum] (0,1.44)--(0.85,2.88);
\draw (0.84,2.88)--(-0.85,2.88) ;

\draw (-1.68,1.44)--(-0.84,0)--(0,1.44);

\draw (0,1.44)--(0.84,0);

\draw (0.84,0)--(1.68,1.44);

\draw (-1.3,2.95) node {$112$};
\draw (1.3, 2.95) node {$113$};
\draw (0,4.6) node   {$111$};
\draw (-2.8,-0.1) node   {$222$};
\draw (2.8,-0.1) node   {$333$};
\draw (2.15,1.5) node   {$133$};
\draw (-2.15,1.5) node   {$122$};
\draw (-0.84,-0.3) node {$223$};
\draw (0.84, -0.3) node {$233$};
\draw (0.45, 1.6) node {$123$};

\end{tikzpicture}
    \caption{The cell $C_{(1,0,1),123}$ in $\cH^3_3$.}
    \label{fig: H33UpTriang}
\end{subfigure}
\hspace{0.02\textwidth}
\begin{subfigure}[b]{0.47\textwidth}
\centering
\begin{tikzpicture}[]

\draw    (0, 4.33) -- (-2.5,0)--(2.5,0)--(0,4.33) ;

\draw    (-1.68, 1.44)--(0,1.44);

\draw[very thick, BlueGreen] (0,1.44)--(1.68, 1.44) ;

\draw  (-0.85,2.88)--(0,1.44)--(0.85,2.88)--(-0.85,2.88) ;

\draw (-1.68,1.44)--(-0.84,0)--(0,1.44);

\draw[very thick, Brown] (0,1.44)--(0.84,0);

\draw[very thick, BurntOrange] (0.84,0)--(1.68,1.44);

\draw (-1.3,2.95) node {$112$};
\draw (1.3, 2.95) node {$113$};
\draw (0,4.6) node   {$111$};
\draw (-2.8,-0.1) node   {$222$};
\draw (2.8,-0.1) node   {$333$};
\draw (2.15,1.5) node   {$133$};
\draw (-2.15,1.5) node   {$122$};
\draw (-0.84,-0.3) node {$223$};
\draw (0.84, -0.3) node {$233$};
\draw (0.45, 1.6) node {$123$};

\end{tikzpicture}
\caption{The cell $C_{(0,0,1),123}$ in $\cH^3_3$.}\label{fig: H33downTriang}
\end{subfigure}
\caption{}\label{fig: mainEx}
\end{figure}

\begin{example}\label{ex: 2D}
Consider the hypersimplicial complex $\cH^3_3$ pictured in Figure \ref{fig: dilatedSimplex}. In Figure \ref{fig: mainEx}, we indicate for each vertex the corresponding element of $S_3(F)$ and color some of the faces for emphasis. The edge with vertices $(1,1,1)$ and $(1,0,2)$ in \textcolor{BlueGreen}{blue-green} in Figure \ref{fig: H33UpTriang} and \ref{fig: H33downTriang} corresponds to the cell $C_{(1,0,1),23}$, which has the following image in the double complex of Figure \ref{fig: doubleCplx}:
$$
    \begin{ytableau}
    \textcolor{BlueGreen}{2} & \textcolor{BlueGreen}{1} &
    \textcolor{BlueGreen}{3}\\
    \textcolor{BlueGreen}{3}
    \end{ytableau} \mapsto
   x_2 \cdot \begin{ytableau}
   1 & 3 & 3\\
   \end{ytableau} - x_3\cdot \begin{ytableau}
   1 & 2 & 3\\
   \end{ytableau}.
$$
The ``up-simplex'' corresponding to the cell $C_{(1,0,1),123}$ colored in Figure \ref{fig: H33UpTriang} has image
$$
    \begin{ytableau}
    1 & 1 & 3\\
    2\\
    3
    \end{ytableau} \mapsto
    x_1 \cdot \begin{ytableau}
    \textcolor{BlueGreen}{2} & \textcolor{BlueGreen}{1} &
    \textcolor{BlueGreen}{3}\\
    \textcolor{BlueGreen}{3}
    \end{ytableau} - 
    x_2 \cdot \begin{ytableau}
    \textcolor{RedOrange}{1} & \textcolor{RedOrange}{1} &
    \textcolor{RedOrange}{3}\\
    \textcolor{RedOrange}{3}
    \end{ytableau} +
    x_3\cdot \begin{ytableau}
    \textcolor{Plum}{1} & \textcolor{Plum}{1} &
    \textcolor{Plum}{3}\\
    \textcolor{Plum}{2}
    \end{ytableau}.
$$
Finally, the ``down-simplex'' of Figure \ref{fig: H33downTriang} corresponding to the cell $C_{(0,0,1),123}$ has image
$$
    \begin{ytableau}
    1 & 3\\
    2\\
    3
    \end{ytableau} \mapsto
\begin{ytableau}
    \textcolor{BlueGreen}{2} & \textcolor{BlueGreen}{1} &
    \textcolor{BlueGreen}{3}\\
    \textcolor{BlueGreen}{3}
    \end{ytableau} - \begin{ytableau}
    \textcolor{Brown}{1} & \textcolor{Brown}{2} &
    \textcolor{Brown}{3}\\
    \textcolor{Brown}{3}
    \end{ytableau} +
    \begin{ytableau}
    \textcolor{BurntOrange}{1} & \textcolor{BurntOrange}{3} &
    \textcolor{BurntOrange}{3}\\
    \textcolor{BurntOrange}{2}
    \end{ytableau}.
$$
In general, the number of elements in the column of a given hook tableaux indicates the homological degree in which it appears in the double complex of Figure \ref{fig: doubleCplx}, or equivalently the dimension (plus one) of the corresponding face in the hypersimplicial complex.
\end{example}

The following simple observation turns out to be critical for applications in Section \ref{sec: dictionary}.

\begin{obs} \label{obs: downTriangStraightening}
All elements of $\bigwedge^2 F \otimes S_{d-1}$ contained in the one-skeleton of a ``down-triangle'' are related by a straightening relation as in Observation \ref{obs:straighteningrelations}.
\end{obs}

\begin{figure}
    \centering
    \hspace{3em}\includegraphics[width = 0.6\textwidth]{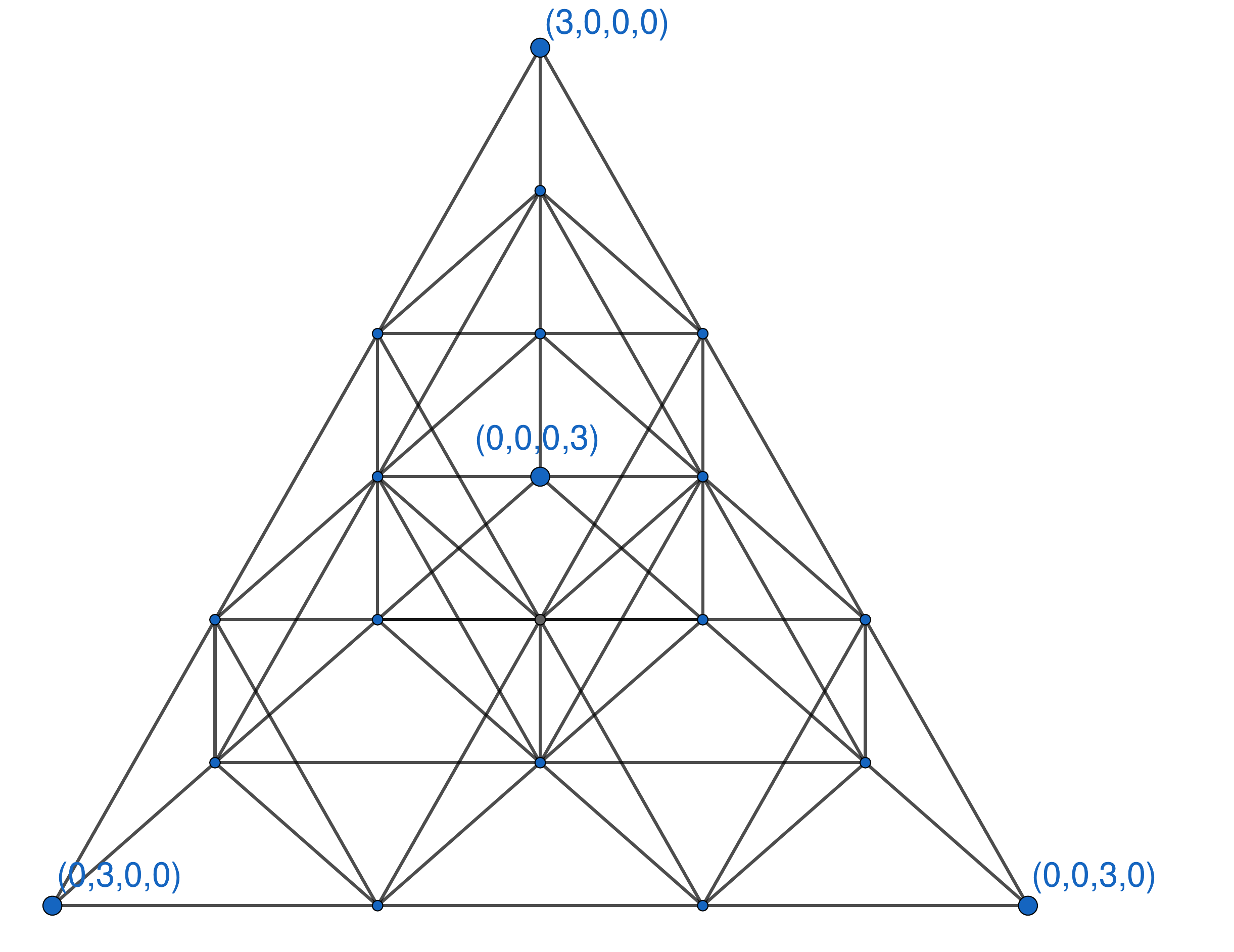}
    \caption{The one-skeleton of $\cH^3_4$.}
    \label{fig:hypersimplex}
\end{figure}
\begin{example}\label{ex: hypersimplexRes}
Figure \ref{fig:hypersimplex} depicts the one-skeleton of $\cH^3_4$. In this case, the hypersimplices which appear as maximal cells are not solely simplices: there are four octahedra in the complex, which correspond to the four possible elements of $\bigwedge^4 F \otimes F$. The image of any one of these octahedra in the double complex of Figure \ref{fig: doubleCplx} is a linear combination of eight tableaux, corresponding to the eight faces of the octahedron: four ``up-triangles'' coming from its image under $\kappa$, and four ``down-triangles'' corresponding to its image under $\kos^\psi\otimes 1$.
\end{example}

\begin{prop}[{see \cite{batzies2002discrete}}]\label{prop: hypersimplexCellRes} Let $\m^d = (x_1,\dots, x_n)^d\subset k[x_1,\dots, x_n]$. Then $\cH^d_n$ defines a multigraded cellular free resolution of $\m^d$.
\end{prop}

Batzies and Welker \cite{batzies2002discrete} use discrete Morse theory to show that the Eliahou--Kervaire resolution for powers of the graded maximal ideal is cellular. We apply their techniques to obtain a minimal cellular resolution isomorphic to the $L$-complex as in Definition \ref{def:Lcomplexes}. To do this, we find an acyclic matching on $\cH^d_n$ distinct from the one in \cite{batzies2002discrete} which has a corresponding Morse complex which supports the $L$-complex. 
\begin{notation}\label{not: firstMin} For any vector $\aa\in\NN^n$, denote by $\first(\aa)$ the first nonzero entry of $\aa$. For an indexing set $J\subset [n]$, denote by $\min J$ the smallest integer appearing in $J$.
\end{notation}

\begin{prop}\label{prop: cellularLcplx} Let $\cH_n^d$ be as in Definition \ref{def: hypersimplex}. Consider the matching $C_{\aa,J} \ra C_{\aa+\epsilon_{\min J}, J\setminus \min{J}}$, where
\begin{enumerate}
    \item $\aa\in\NN^n$,
    \item $J\subset [n]$ is such that $2\leq d-\abs{\aa} \leq \abs{J}-1$, and
    \item $\min J \leq \first(\aa)$.
\end{enumerate}
Then this is an acyclic homogeneous matching $A$ on $\cH^d_n$ as in Construction \ref{const: DMT}. 
\end{prop}

\begin{proof}
$A$ is a matching because cells $C_{\aa,J}$ on the left hand side must satisfy $\min J \leq \first(\aa)$ while the cells $C_{\aa',J'}$ on the right hand side must satisfy $\min J' > \first(\aa')$. Suppose, seeking contradiction, that $A$ contains a cycle. Observe that $\mdeg(C_{\aa,J})$ must be weakly decreasing along every directed edge of the cycle, so in particular it must be constant. Observe that every element $C_{\aa,J}$ at the head of an edge directed upwards in $E^A$ satisfies $\min J \leq \first(\aa)$, but then every element at the head of an arrow pointing ``down'' from one of these elements must have the same $\min J$. But this element $C_{\aa',J'}$ must also be at the tail of some other element of $E^A$ pointing upwards, so it must also satisfy that $\min J' < \first(\aa)$, which is a contradiction.
\end{proof}

By characterizing the cells untouched by the acyclic matching in the previous proposition, we obtain the following corollary.

\begin{cor}\label{cor: criticalCells}
Let $\tilde \cH^d_n = (\cH^d_n)_A$ denote the corresponding Morse complex after applying the acylic matching from Proposiiton \ref{prop: cellularLcplx}. Then the $A$-critical cells of $\cH^d_n$ are:
\begin{enumerate}
    \item the $0$-cells $C_{\aa,\emptyset}$, where $\aa\in \NN^n\cap\Delta (n,d)$, and
    \item the cells $C_{\aa,J}$ such that $\min J \leq \first(\aa)$ and $\abs{\aa} = d-1$.
\end{enumerate}
\end{cor}

We conclude this section with our main result, which states that the $L$-complex from Definition \ref{def:Lcomplexes} is supported on a CW-complex.

\begin{prop}\label{prop: LcplxCellular} Let $\tilde\partial$ denote the differential of the Morse complex $\tilde H_n^d = (H_n^d)_A$ and $\tilde C_{\aa, J}$ the cell in $\tilde H_n^d$ corresponding to the $A$-critical cell $C_{\aa,J}$ of $\cH^d_n$. Then 
$\tilde \cH_n^d$ supports a minimal linear cellular resolution of a power of the graded maximal ideal which is isomorphic to the $L$-complex from Definition \ref{def:Lcomplexes}.
\end{prop}

\begin{proof}
The $A$-critical cells of $\cH^d_n$ are exactly the $0$-cells $C_{\aa,\emptyset}$ for $\aa\in \Delta_n(d)\cap \NN^n$ and all cells $C_{\aa,J}$ such that $\min J\leq \min \aa$ and $\abs{\aa} = d-1$; in particular, the critical cells correspond to exactly those standard hook tableaux which are basis elements of the modules in the $L$-complex. Let $C_{\aa,J}$ be a critical cell with $p = \min J$ and $q = \min\{i\mid i\in \supp(\aa)\}$. Then the differential $\partial$ of $\cH_n^d$ applied to $C_{\aa,J}$ has at most one nonstandard tableau in its image, which would be $C_{\aa, J_p}$. This element is matched with $C_{\aa-\epsilon_q, J_p\cup q}$, which has an image under $\partial$ consisting of standard tableaux of the same shape and multidegree as $C_{\aa, J_p}$, and potentially some tableaux corresponding to elements which are matched with elements one dimension lower and therefore do not appear in $\tilde \partial$. In particular, after homogenization, $\tilde \partial$ corresponds exactly to the differential of the $L$-complex in Definition \ref{def:Lcomplexes}.
\end{proof}
\section{Polarizations of Powers of Graded Maximal Ideals}\label{sec: pols}

The material in this section is a summary of the combinatorial characterization of polarizations of powers of the graded maximal ideal $\m$ in a polynomial ring given by Almousa, Fl\o ystad, and Lohne in \cite{polarizations}. The idea is to put a set of partial orders $\geq_i$ on the generators of $\mathfrak{m}^d$ and view a ``potential polarization'' as a set of isotone maps. One may visualize these potential polarizations as a graph of linear syzygies among the generators of the ``potentially polarized'' ideal, which is in fact a subgraph of the one-skeleton of the hypersimplicial complex $\cH_n^d$ from the previous section. Theorem \ref{thm: polsSyzEdges} gives a complete characterization of which of these isotone maps give ``honest'' polarizations in terms of their graphs of linear syzygies. The main result of this section is given by Theorem \ref{thm: polsSyzEdges}, which will be retranslated in the context of Schur modules in Section \ref{sec: dictionary}.

Much of the notation introduced in this section will be used in later sections without reference. We begin this section with the definition of a polarization. 
Intuitively, polarizations can be used to replace any monomial ideal with a \emph{squarefree} monomial that is homologically indistinguishable.
One of the main insights of \cite{polarizations} is the fact that there are many ways to polarize an ideal, and the family of all polarizations of a given monomial ideal can be highly nontrivial. However, in the case of the ideal $\m^d$, a full characterization of all such polarizations is possible; we recall this result in Proposition \ref{thm: polsSyzEdges}.

\begin{definition}[Polarization]\label{def: polarization} Let $I$ be a monomial ideal in a polynomial ring $S = k[x_1,\dots, x_n]$ over a field $k$, and let $d_i$ be the highest power of $x_i$ that appears in a minimal generator of $I$. Set $\Xv_i = \{x_{i1},\dots, x_{i d_i}\}$ for all $i$, and define the polynomial ring $\tilde S = k[\Xv_1,\dots, \Xv_n]$ in the union of all these variables. Then $\tilde I\subset \tilde S$ is a \emph{polarization} of $I$ if
$$
\sigma = \bigcup_{i=1}^n \{x_{i1}-x_{i2}, x_{i1} - x_{i3}, \dots,x_{i1}-x_{i d_i} \}
$$
is a regular $\tilde S / \tilde I$-sequence and $\tilde I \otimes \tilde S/\langle \sigma \rangle \cong I$.
\end{definition}

The first incarnation of polarizations appeared in Hartshorne's thesis (see \cite{hartshorne1966connectedness}), where he used what he called ``distractions'' to prove the connectedness of the Hilbert scheme. 

\begin{setup}\label{set: pols}
Fix integers $n$ and $d$, and let $S = k[x_1,\dots, x_n]$ be a polynomial ring over a field $k$. Let $\Xv_i = \{x_{i1},\dots, x_{i d}\}$ be a set of variables, and let $\tilde S = k[\Xv_1,\dots, \Xv_n]$ be a polynomial ring in the union of all these variables. Denote by $\mathfrak{m} = (x_1,\dots, x_n)$ the graded maximal ideal of $S$.

Let $\Delta(n,d)$ be the dilated $(n-1)$-simplex from Definition \ref{not: dilatedSimplex}. Denote by $\Delta^\ZZ(n,d) = \Delta(n,d)\cap \ZZ^n$ the set of lattice points of the dilated simplex $d\cot \Delta_{n-1}$, i.e., the set of tuples $\aa = (a_1,\dots, a_n)$ of non-negative integers with $\sum_i^n a_i = d$. Denote by $\sk_1(\cH^d_n)$ the one-skeleton of the hypersimplicial complex $\cH^d_n$ from Definition \ref{def: hypersimplex}.
\end{setup}

\begin{remark}
Observe that the elements of $\Delta^\ZZ(n,d)$ are exactly the exponent vectors of the minimal generating set of the ideal $\mathfrak{m}^d$.
\end{remark}

\begin{notation} Let $\epsilon_i\in\NN^n$ be the $i$th unit vector in $\NN^n$.
For a given $\aa$, denote by $\Supp(\aa)$ the \emph{support} of $\aa$, that is, the set of all $i$ such that $a_i>0$. 
If $B$ is a subset of $[n]$, denote by $\one_B$ the $n$-tuple $\sum_{i\in B} \epsilon_i$. For example, if $B = [n]$, then $\one_B = (1,\dots, 1)$.
\end{notation}

\begin{definition}\label{def: partialOrder} Adopt notation and hypotheses of Setup \ref{set: pols}. Fix an index $1\leq i\leq n$. Define $(\Delta^\ZZ(n,d),\geq_i)$ to be the poset with ground set $\Delta^\ZZ(n,d)$ and partial order $\geq_i$ such that $\bb\geq_i \aa$ if $b_i\geq a_i$ and $b_j\leq a_j$ for $j\neq i$. 
\end{definition}

\begin{obs}\label{obs: gradedPoset}
The partial order $\geq_i$ as in Definition \ref{def: partialOrder} is graded, where $\aa\in\Delta^\ZZ(n,d)$ has rank $a_i$.
\end{obs} 

In the following definitions, we introduce some key subgraphs of the one-skeleton of $\cH^d_n$ which will be critical for the combinatorial characterization of polarizations of $\m^d$. Complete down-graphs will also be significant players in Section \ref{sec: dictionary}.

\begin{definition}[Complete down-graph]\label{def: completeDownGraph} Given $\cc\in\Delta^\ZZ(n,d+1)$ and $i,j\in\Supp(\cc)$, there is an edge between $\cc-\epsilon_i$ and $\cc-\epsilon_j$ in $\sk_1(\cH_n^d)$ denoted $(\cc;i,j)$. Every edge in $\sk_1(\cH_n^d)$ can be realized as an edge $(\cc;i,j)$ for unique $\cc,i,$ and $j$. An $n$-tuple $\cc\in\Delta^\ZZ(n,d+1)$ induces a subgraph of $\sk_1(\cH_n^d)$ called the \emph{complete down-graph} $D(\cc)$ on the points $\cc-\epsilon_i$ for $i\in\Supp(\cc)$.
If $R\subseteq [n]$, denote by $D_R(\cc)$ the complete graph with edges $(\cc;r,s)$ for $r,s\in R$.
\end{definition}

\begin{definition}[Complete up-graph]\label{def: completeUpGraph} Any $\aa\in\Delta^\ZZ(n,d-1)$ also determines a subgraph of $\sk_1(\cH_n^d)$: the \emph{complete up-graph} $U(\aa)$ consisting of points $\aa+\epsilon_i$ for $i=1,\dots, n$ with edges $(\aa+\epsilon_i+\epsilon_j;i,j)$ for $i\neq j$.
\end{definition}

\begin{remark} The complete down-graph $D(\cc)$ induces a simplex of full dimension $d-1$ if and only if $c_i\geq 1$ for all $i$, i.e., $\cc$ has full support. For each $\aa$ in $\Delta^\ZZ(n,d-1)$, the induced simplex of the up-graph $U(\aa)$ always has full dimension $d-1$.
\end{remark}

\begin{example}\label{ex: Delta_3(3)} The one-skeleton of $\cH_3^3$ pictured in Figure \ref{fig: dilatedSimplex} has three ``complete down-triangles'' with full support corresponding to the vectors $(2,1,1), (1,2,1)$, and $(1,1,2)$ in $\Delta^\ZZ(n,d+1)$. It also has six ``complete up-triangles''.
\end{example}

The maps in the following construction will be play an important role in the combinatorial characterization of $\m^d$.

\begin{construction}\label{const: pols}
Adopt notation and hypotheses of Setup \ref{set: pols}. Let $\cB_d$ be the Boolean poset on $[d]$ and $\{X_i\}_{1\leq i\leq n}$ be a set of rank-preserving isotone maps
\begin{align*}
    X_i: (\Delta^\ZZ(n,d),\leq_i) \ra \cB_d.
\end{align*}
For any $\aa\in\Delta^\ZZ(n,d)$, let $m_i(\aa) = \prod_{j\in X_i(\aa)} x_{ij}$ and $m(\aa) = \prod_{i=1}^n m_i(\aa)$. Let $J$ be the ideal in $k[\Xv_1,\dots, \Xv_n]$ generated by the $m(\aa)$.
\end{construction}

\begin{definition}[Linear syzygy edge]\label{def: linSyzEdge} 
Let $(\cc;i,j)$ be an edge of $\sk_1(\cH_n^d)$, where $\cc\in\Delta^\ZZ(n,d+1)$. Then $(\cc ; i,j)$ is a \emph{linear syzygy edge} (or \emph{LS-edge}) if there is a monomial $\mm$ of degree $d-1$ such that
$$
m(\cc-\epsilon_i) = x_{jr}\cdot \mm \quad \textrm{and} \quad m(\cc-\epsilon_j) = x_{is}\cdot \mm,
$$
for suitable variables $x_{jr}\in \Xv_j$ and $x_{is}\in\Xv_i$. This edge gives a linear syzygy between the monomials $m(\cc-\epsilon_i)$ and $m(\cc-\epsilon_j)$. Equivalently, in terms of the isotone maps,
$$
X_p(\cc-\epsilon_i) = X_p(\cc-\epsilon_j)
$$
for every $p\neq i,j$. Observe that both $m_i(\cc-\epsilon_i)$ and $m_j(\cc-\epsilon_j)$ are common factors of $m(\cc-\epsilon_i)$ and $m(\cc-\epsilon_j)$.
\end{definition}

\begin{notation}\label{not: LS}
For any $\cc\in\Delta^\ZZ(n,d+1)$, let $\LS(\cc)$ be the set of linear syzygy edges in the complete down-graph $D(\cc)$. For any $\bb\in\Delta^\ZZ(n,d-1)$, denote by $\LS_i(\bb)$ the set of linear syzygy edges in the induced subgraph of $U(\bb)$ on the vertex set $\{\bb+\epsilon_j \mid j\neq i\}$.
\end{notation}

Sometimes, one may wish to consider whether two elements of $\Delta^\ZZ(n,d)$ would share a linear syzygy edge with respect to a \textit{subset} of $[n]$.

\begin{definition}[$R$-linear syzygy edge]\label{def: R-LSedges} Let $R\subseteq [n]$ and $\cc\in\Delta^\ZZ(n,d+1)$ with $R$ contained in the support of $\cc$. Let $r,s\in R$. Define $(\cc;r,s)$ to be an \emph{$R$-linear syzygy edge} if
$$
X_p(\cc-\epsilon_r) = X_p(\cc-\epsilon_s) \text{ for } p\in R\setminus \{r,s\}.
$$
By the isotonicity of the $X_p$, for $p = r,s$,
$$
X_r(\cc - \epsilon_r)\subseteq X_r(\cc - \epsilon_s),\quad X_s(\cc-\epsilon_s)\subseteq X_s(\cc-\epsilon_r).
$$
Let $D_R(\cc)$ be the complete graph with edges $(\cc;r,s)$ for
$r,s \in R$.
\end{definition}

The following lemma will be particularly useful in Section \ref{sec: dictionary}.

\begin{lemma}\label{lem: R-LSedge} Let $\cc\in\Delta^\ZZ(n,d+1)$. If the set of linear syzygy edges in $\LS(\cc)$ contains a spanning tree for $D(\cc)$, then for each $R\subseteq \supp(\cc)$, the set of $R$-linear syzygy edges contains a spanning tree for $D_R(\cc)$.
\end{lemma}

We conclude this section by presenting the main theorem of \cite{polarizations}: a complete combinatorial characterization of all polarizations of $\m^d$ in terms of their graphs of linear syzygies. 

\begin{theorem}[\cite{polarizations}]\label{thm: polsSyzEdges} Adopt notation and hypotheses of Setup \ref{set: pols} and Construction \ref{const: pols}. A set of isotone maps $X_1,\dots, X_n$ determines a polarization of the ideal $(x_1,\dots, x_n)^d$ if and only if for every $\cc\in\Delta^\ZZ(n,d+1)$, the linear syzygy edges $\LS(\cc)$ contain a spanning tree for the down-graph $D(\cc)$.
\end{theorem}

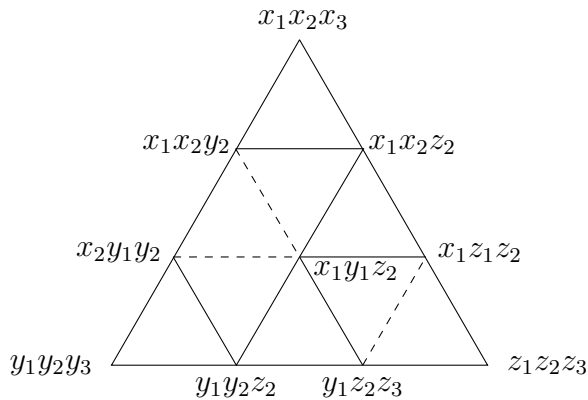
\begin{figure}
\begin{tikzpicture}

\draw  (0, 4.33) -- (-2.5,0)--(2.5,0)--(0,4.33) ;

\draw[dashed]  (-1.68, 1.44)--(0,1.44);

\draw (0,1.44)--(1.68, 1.44) ;

\draw[dashed] (-0.85,2.88)--(0,1.44);

\draw (0,1.44)--(0.85,2.88)--(-0.85,2.88) ;

\draw (-1.68,1.44)--(-0.84,0)--(0,1.44);

\draw (0,1.44)--(0.84,0);

\draw[dashed] (0.84,0)--(1.68,1.44);

\draw (-1.5,2.93) node {$x_1x_2y_2$};
\draw (1.5, 2.93) node {$x_1x_2z_2$};
\draw (0.05,4.6) node   {$x_1x_2x_3$};
\draw (-3.3,0) node   {$y_1y_2y_3$};
\draw (3.3,0) node   {$z_1z_2z_3$};
\draw (2.4,1.5) node   {$x_1 z_1z_2$};
\draw (-2.4,1.5) node   {$x_2y_1y_2$};
\draw (-0.84,-0.3) node {$y_1y_2z_2$};
\draw (0.84, -0.3) node {$y_1z_2z_3$};
\draw (0.75, 1.25) node {$x_1y_1z_2$};

\end{tikzpicture}
\caption{An example of a polarization of $(x,y,z)^3$.}
\label{fig: LS-xyz}
\end{figure}

\begin{example}\label{ex: pol} Figure \ref{fig: LS-xyz} depicts the graph of linear syzygies for a polarization of $(x,y,z)^3$. Notice that at most one edge is removed from each down-triangle, so it satisfies the spanning tree condition of Theorem \ref{thm: polsSyzEdges}.
\end{example}
\section{Hook Tableaux and Polarizations}\label{sec: dictionary}

The goal of this section is to provide a dictionary between the notation and terminology introduced in Section \ref{sec: pols} and the Schur modules appearing in the $L$-complexes of Section \ref{sec: Lcomplex}. More precisely, we give a new combinatorial characterization of all polarizations of powers of the graded maximal ideal in terms of generating sets of the Schur module $L^1_d(F)$. The main result of this section is Theorem \ref{thm: spanningTree}, which shows that the spanning tree condition of Theorem \ref{thm: polsSyzEdges} is equivalently asking that the Young tableaux canonically associated to the linear syzygy edges form a basis for the associated Schur module. 

The actual dictionary for translating between the different aforementioned frameworks is given by Proposition \ref{prop: dictionary}; these results will be employed in Section \ref{sec: boundedGens} to extend the results of \ref{def: polarization} to the case of restricted powers. 

\begin{setup}\label{set: dictionary}
Fix integers $n$ and $d$, and let $S = k[x_1,\dots, x_n]$ be a polynomial ring over a field $k$. Let $\Xv_i = \{x_{i1},\dots, x_{i d}\}$ be a set of variables, and let $\tilde S = k[\Xv_1,\dots, \Xv_n]$ be a polynomial ring in the union of all these variables. Denote by $\mathfrak{m} = (x_1,\dots, x_n)$ the graded maximal ideal of $S$.

Let $\Delta(n,d)$ be the dilated $(n-1)$-simplex from Definition \ref{not: dilatedSimplex}. Denote by $\Delta^\ZZ(n,d) = \Delta(n,d)\cap \ZZ^n$ the set of lattice points of the dilated simplex $d\Delta_{n-1}$, i.e., the set of tuples $\aa = (a_1,\dots, a_n)$ of non-negative integers with $\sum_i^n a_i = d$. Denote by $\sk_1(\cH^d_n)$ the one-skeleton of the hypersimplicial complex $\cH^d_n$ from Definition \ref{def: hypersimplex}.

Let $L^a_b(F)$ be the Schur module defined in Setup \ref{set:Lcomplexsetup}.
\end{setup}

\begin{prop}\label{prop: dictionary} Adopt notation and hypotheses of Setup \ref{set: dictionary}. Then:
\begin{enumerate}[label = \textup{(}\alph*\textup{)}]
    \item There exists a bijection $\psi_{n,d}$ from $\Delta^\ZZ(n,d)$ to $S_d(F)$.
    \item For any pair $\cc\in\Delta^\ZZ(n,d+1)$ and $R\subseteq \Supp(\cc)$ such that $\abs{R} = t$, the complete subgraph $D_R(\cc)$ (see Definition \ref{def: R-LSedges}) corresponds to a unique element of $\bigwedge^t F \otimes S_{d-t+1}(F)$.
    \item There exists a bijection $\theta_{n,d}$ from the set of edges of $\Delta^\ZZ(n,d)$ to $\bigwedge^2 F\otimes S_{d-1}(F)$.
\end{enumerate}
\end{prop}

\begin{proof}
For (a), the map
\begingroup\allowdisplaybreaks
\begin{equation}
    \psi_{n,d}: \Delta^\ZZ(n,d) \ra S_d(F)
\end{equation}
\endgroup
such that $\psi_{n,d}(\aa) = f^\aa$
gives the desired bijection.

For (b), Let $\cc\in\Delta^\ZZ(n,d+1)$. If $R\subseteq \Supp(\cc)$, define
\begin{equation}
\omega_{n,d}(\cc,R) \coloneqq f_R \otimes \psi_{n,d-\abs{R}+1}(\cc_R)
\end{equation}
where $\cc_R = \cc - \sum_{i\in R} \epsilon_i$. In particular, if $R = \Supp(\cc)$, then the map
\begin{equation}
\omega_{n,d}: \Delta^\ZZ(n,d+1) \ra \bigwedge^t F \otimes S_{d-t+1}(F) \\
\end{equation}
such that $\omega_{n,d}(\cc) = f_{\Supp (\cc)} \otimes \psi_{n,d-t+1}(\cc')$ is a bijection between the down-triangles of $\sk_1(\cH_n^d)$ and $\bigwedge^t F \otimes S_{d-t+1}(F)$, where $\cc' = \cc - \sum_{i\in \Supp(\cc)} \epsilon_i$ and $t = \abs{\Supp(\cc)}$.

For (c), If $(\cc;i,j)$ is an edge in $\sk_1(\cH_n^d)$, then the map
\begin{equation}
\theta_{n,d}(\cc; i,j) \coloneqq f_i\wedge f_j \otimes \psi_{n,d-1}(\cc- \epsilon_i-\epsilon_j)
\end{equation}
gives a bijection between edges of $\sk_1(\cH_n^d)$ and elements of $\bigwedge^2 F \otimes S_{d-1}(F)$.
\end{proof}

\begin{chunk} Let $\cc\in\Delta^\ZZ(n,d+1)$ be a down-triangle and let $R \subseteq \Supp(\cc)$. Set $T = \omega(\cc,R)\in \bigwedge^{\abs{R}} F \otimes S_{d-\abs{R}-1}(F)$.
The tableaux appearing in the image $\kappa_{\abs{R},n-\abs{R}+1}(T)$ (where $\kappa_{a,b}$ is defined in Setup \ref{set:Lcomplexsetup}) correspond exactly to the down-triangles $D_P(\cc)$ such that $\abs{P} = \abs{R}-1$. In particular, as noted in Observation \ref{obs: downTriangStraightening}, if $\abs{R} = 3$, $\kappa_{3,n-2}(T)$ corresponds to a linear combination of the labels appearing on the three edges of $D_R(\cc)$.
\end{chunk}

\begin{definition}\label{def: tabsChi} Let $\chi = \{X_i\}$ denote a set of isotone maps
$$
X_i : (\Delta^\ZZ(n,d), \geq_i) \ra \cB_d
$$
as in Construction \ref{const: pols}. Let $\LS_\chi$ be the set of linear syzygy edges after applying $\chi$ to the generators of $\m^d$ as in Definition \ref{def: linSyzEdge}.
Denote by $\tabs(\chi)$ the set of tableaux in $\bigwedge^2 F \otimes S_{d-1}(F)$ associated to the edges $\LS_\chi$ via the correspondence in \ref{chunk: hypersimplexHooks}.
\end{definition}

\begin{theorem}\label{thm: spanningTree} Adopt notation and hypotheses of Setup \ref{set: dictionary}. Let $\chi = \{X_i\}$ denote a set of isotone maps as in Construction \ref{const: pols}.
Then the following are equivalent:
\begin{enumerate}
    \item The elements of $\tabs(\chi)$ span the module $L_d^1(F)$.
    \item For every $\cc\in\Delta^\ZZ(n,d+1)$, $\LS(\cc)$ contains a spanning tree of the complete down-graph $D(\cc)$.
    \item The set of isotone maps $X_1,\dots, X_n$ determine a polarization of $(x_1,\dots, x_n)^d$.
\end{enumerate}
\end{theorem}

\begin{proof}
Note that (2)$\iff$(3) is Theorem \ref{thm: polsSyzEdges}.

\noindent\textbf{(2)$\implies$(1):} 
Take $\cc\in\Delta^\ZZ(n,d+1)$ and let $\aa = \cc-\epsilon_i$ and $\bb = \cc-\epsilon_j$ for some $i,j\in\Supp(\cc)$.
If $LS(\cc)$ contains a spanning tree, then for any two vertices $\aa$ and $\bb$ in $\cc$, there exists a path in $LS(\cc)$ connecting them.
It suffices to show that for any $\aa$ and $\bb$ in $D(\cc)$, $\theta_{n,d}(\cc;i,j)$ is in the span of $\tabs(\chi)$.

Proceed by induction on $k$, the number of edges in the shortest path from $\aa$ to $\bb$. If $k=1$, then the tableau corresponding to the edge between $\aa$ and $\bb$ is in $\tabs(\chi)$. Now assume that for any two vertices $\cc-\epsilon_i$ and $\cc-\epsilon_j$ such that the shortest path in $LS(\cc)$ between them is length $k$, the tableau $\theta_{n,d}(\cc;i,j)$ is a linear combination of elements in $\tabs(\chi)$. Let $\aa$ and $\bb$ be two vertices of $D(\cc)$ such that the shortest path between them is length $k+1$, i.e., there is a set of vertices
$$
\aa = \dd^1,\dots, \dd^{k+2} = \bb
$$
such that each $\dd^j$ is a vertex in $D(\cc)$ and each pair $(\dd^j,\dd^{j+1})$ is connected by an edge in $\LS(\cc)$. The length of the shortest path between $\aa$ and $\dd^{k+1}$ is $k$, so by the induction hypothesis, the tableau labeling the edge between them is spanned by the elements of $\tabs(\chi)$. If $\aa = \cc-\epsilon_{i_1},\bb = \cc-\epsilon_{i_2},$ and $\dd^{k+1} = \cc-\epsilon_{i_3}$, set $R = \{i_1, i_2, i_3\}$ and consider the ``smaller'' down-triangle $D_R(\cc)$ (see Definition \ref{def: R-LSedges}). By Observation \ref{obs: downTriangStraightening}, the tableau corresponding to the edge between $\aa$ and $\bb$ is a linear combination of the tableaux corresponding to the other two edges of $D_R(\cc)$, which in turn have been shown to be linear combinations of elements of $\tabs(\chi)$, hence proving the claim.

\noindent\textbf{(1)$\implies$(2):} 
Let $\cc\in\Delta^\ZZ(n,d+1)$ be a complete down-graph in $\sk_1(\cH_n^d)$. It suffices to show that for any two vertices $\aa = \cc-\epsilon_i$ and $\bb = \cc-\epsilon_j$, there exists a path from $\aa$ to $\bb$ by edges labeled by tableaux in $\tabs(\chi)$. Since $\tabs(\chi)$ must contain a basis of the module $L^1_d(F)$, one may assume that $\tabs(\chi)$ is a basis itself.
Suppose $\theta_{n,d}(\cc;i,j)$ is not in $\tabs(\chi)$. The result follows from the following claims:
\begin{enumerate}[(i)]
    \item There exists some $R = \{i,j,\ell_1\}\subseteq \Supp(\cc)$ such that $D_R(\cc)$ has at least one edge labeled by an element of $\tabs(\chi)$.
    \item Suppose $(\cc;i,\ell_1)$ is the unique edge of $D_R(\cc)$ labeled by an element of $\tabs(\chi)$. Then there exists some $P = \{\ell_1, \ell_2, j\}\subseteq \Supp(\cc)$ such that at least one edge of $D_P(\cc)$ has a label appearing in $\tabs(\chi)$ and $\ell_2\neq i$.
    \item Let $\Gamma(\cc)$ be the subgraph of $D(\cc)$ with edges labeled by elements of $\tabs(\chi)$. If $\Gamma(\cc)$ contains a cycle, then it corresponds to a linearly dependent subset of $\tabs(\chi)$.
\end{enumerate}

To see (i), suppose no tableaux corresponding to edges in any possible $D_R(\cc)$ are in $\tabs(\chi)$. Then no tableaux in any of the possible straightening relations containing $\theta_{n,d}(\cc;i,j)$ coming from the image of any of the $\omega_{n,d}(\cc,R)$ under $\kappa_{3,n-2}$ appear in $\tabs(\chi)$. Hence, $\theta_{n,d}(\cc;i,j)$ is not in the span of $\tabs(\chi)$.

For (ii), observe that the image of $\omega_{n,d}(\cc,R)$ under $\kappa_{3,n-2}$ gives that $\theta_{n,d}(\cc;i,j)$ is in the span of $\theta_{n,d}(\cc; i,\ell_1)$ and $\theta_{n,d}(\cc; j,\ell_1)$. By assumption, $\theta_{n,d}(\cc;i,\ell_1)\in\tabs(\chi)$. Suppose every other $D_P(\cc)$ with $\abs{P} = 3$ containing $\theta_{n,d}(\cc;j,\ell_1)$ has no edges labeled by tableaux in $\tabs(\chi)$. Then $\theta_{n,d}(\cc; j,\ell_1)$ is not in the span of $\tabs(\chi)$, contradicting the assumption that $\tabs(\chi)$ spans $L^1_d(F)$.

To check (iii), proceed by induction on the length $k$ of the cycle. In the base case where $k=3$, the cycle forms the edges of a down-triangle $D_R(\cc)$ where $\abs{R}=3$. The three tableaux labeling the edges of $D_R(\cc)$ make up the straightening relation from the image of $\omega_{n,d}(\cc,R)$ under $\kappa_{3,n-2}$, implying they are linearly dependent. This contradicts the assumption that $\tabs(\chi)$ forms a basis. Now assume that any cycle of length $k$ in $\Gamma(\cc)$ gives a linearly dependent subset of $\tabs(\chi)$, and suppose there is a cycle $\aa^1,\dots, \aa^{k+1}$ such that the edges $(\aa^j,\aa^{j+1})$ and $(\aa^1,\aa^k)$ are labeled by tableaux in $\tabs(\chi)$ and each $\aa^j = \cc - \epsilon_{i_j}$. Let $R = \{i_1, i_k, i_{k+1}\}$. By the induction hypothesis, $\theta_{n,d}(\cc; i_1, i_k)$ is equal to a linear combination of tableaux $\theta_{n,d}(\cc; i_k, i_{k+1})$ and $\theta_{n,d}(\cc; i_1, i_{k+1})$; but also by the induction hypotheses, it is a linear combination of tableaux of the form $\theta_{n,d}(\cc; i_j, i_{j+1})$ where $1\leq j\leq k$. Therefore, the cycle corresponds to a linearly dependent subset of $\tabs(\chi)$.

With claims (i)-(iii) established, iterate the following process.
Choose a triangle $D_R(\cc)$ such that $R = \{i,j,\ell\}$ and some edge of $D_R(\cc)$ is labeled by an element of $\tabs(\chi)$. If both $\theta_{n,d}(\cc;i,\ell_1)$ and $\theta_{n,d}(\cc;\ell_1,j)$ are in $\tabs(\chi)$, the claim follows. Suppose that only one of these tableaux are in $\tabs(\chi)$; without loss of generality, assume that $\theta_{n,d}(\cc; i,\ell)\in \tabs(\chi)$. By (ii), there is some triangle $D_P(\cc)$ such that $P = \{\ell_1, \ell_2, j\}$ and at least one edge of $D_P(\cc)$ has a label appearing in $\tabs(\chi)$. If both $\theta_{n,d}(\cc;\ell_1,\ell_2)$ and $\theta_{n,d}(\cc; \ell_2, j)$ appear in $\tabs(\chi)$, then the claim follows. Otherwise, repeat this process. This process must terminate by (iii), giving the desired path.
\end{proof}

\section{Cellular resolutions and Polarizations of Restricted Powers of the Graded Maximal Ideal}\label{sec: boundedGens}

In this section, we extend the results from Sections \ref{sec: pols} and \ref{sec: dictionary} to the case of so-called \emph{restricted powers} of the graded maximal ideal. This class of ideals comes from bounding the multidegrees appearing in the generators of $\m^d$. We show in Proposition \ref{prop: boundedPowersCellRes} that these ideals also have a minimal, linear, cellular resolution arising as a subcomplex of the $L$-complex. In addition, we give two combinatorial characterizations of polarizations of this class of ideals: one in terms of their graphs of linear syzygies, and one in terms of spanning sets of a submodule of the associated Schur module. These characterizations are given by Theorem \ref{thm: generalizedPols}.

\begin{setup}\label{set: boundedGens} Let $I$ be a monomial ideal in a polynomial ring $S = k[x_1,\dots, x_n]$ over a field $k$ with generating set $\cG(I) = (m_1,\dots, m_r)$. Let $\uu\in I$ be a monomial, and define the monomial ideal $I_{\leq \uu}$ to be the ideal generated by monomials $\{m_i\mid m_i \text{ divides } \uu\}$.

Let $d_i$ be the highest power of $x_i$ that appears in a minimal generator of $I$. Set $\Xv_i = \{x_{i1},\dots, x_{i d_i}\}$ for all $i$, and define the polynomial ring $\widetilde{S} = k[\Xv_1,\dots, \Xv_n]$ in the union of all these variables. Observe that $\widetilde S$ has a $\ZZ^n$-grading induced by the first indices of the variables in $\widetilde S$.

Let $\widetilde I$ be a polarization of $I$ as in Definition \ref{def: polarization}. Define $\widetilde I_{\leq \uu}$ to be generated by those elements of $\widetilde I$ with $\ZZ^n$-multidegree bounded above by $\uu$.
\end{setup}

The following useful proposition was observed in \cite{gasharov2002resolutions}.

\begin{prop}\label{prop: subresolution} Let $I$ be a monomial ideal in a polynomial ring $S$. Fix a multihomogeneous basis of a multigraded free resolution $\FF_I$ of $S/I$. Denote by $\FF_I(\leq \uu)$ the subcomplex of $\FF$ generated by multihomogeneous basis elements of multidegrees dividing $\uu$.
\begin{enumerate}[(1)]
    \item The subcomplex $\FF_{I}(\leq \uu)$ is a multigraded free resolution of $S/(I_{\leq \uu})$. 
    \item If $\FF_I$ is a minimal multigraded free resolution of $S/I$, then $\FF_I(\leq \uu)$ is independent of the choice of basis.
    \item If $\FF_I$ is a minimal multigraded free resolution of $S/I$, then the resolution $\FF_I(\leq \uu)$ is also minimal.
\end{enumerate}
\end{prop}

\begin{prop}\label{prop: restrictedPols} Adopt Setup \ref{set: boundedGens}. Then $\widetilde I_{\leq \uu}$ is a polarization of $I_{\leq \uu}$.
\end{prop}

\begin{proof} Let $\widetilde \FF$ be a free resolution of $\widetilde I$. Then $\widetilde \FF$ is a multigraded resolution with respect to the $\ZZ^n$-multigrading where each variable $x_{i,j}$ in $\Xv_i$ has multidegree $e_i$. By the definition of a polarization (Definition \ref{def: polarization}), $\widetilde \FF\otimes \widetilde S/\sigma \cong \FF$, where $\sigma$ is a regular $\widetilde S / \widetilde I$-sequence and $\FF$ is a minimal free resolution of $I$. By Proposition \ref{prop: subresolution}, both $\FF(\leq \uu)$ and $\widetilde \FF(\leq \uu)$  are minimal multigraded free resolutions of $I_{\leq \uu}$ and $\widetilde I_{\leq\uu}$, respectively; in particular, one has that
\begin{equation}\label{eq: exactnessPol}
\widetilde \FF(\leq \uu) \otimes \widetilde S/\sigma \cong \FF(\leq \uu).
\end{equation}
It remains to check that $\sigma$ is indeed a regular sequence on $\widetilde S / \widetilde{I}_{\leq \uu}$. This follows immediately from the string of isomorphisms:
$$H_\bullet (\widetilde \FF(\leq \uu) \otimes \widetilde S/\sigma ) \cong \tor_\bullet^{\widetilde S} (\widetilde S / \widetilde{I}_{\leq \uu} , \widetilde{S} / \sigma) \cong H_\bullet(\widetilde S / \widetilde{I}_{\leq \uu} \otimes K (\sigma)_\bullet ),$$
where $K (\sigma)_\bullet$ denotes the Koszul complex on $\sigma$. Since $\widetilde F(\leq \uu) \otimes \widetilde S/\sigma \cong \FF(\leq \uu)$ is acyclic, it follows that $H_{>0} (\widetilde{S} / \widetilde{I} \otimes K (\sigma)_\bullet ) = 0$, so $\sigma$ is regular on $\widetilde{S} / \widetilde{I}$ by, for instance, \cite[Theorem 14.7]{peeva2010graded}.
\end{proof}

We apply these results to extend results on cellular resolutions and polarizations from powers of the maximal ideal to so-called \emph{restricted powers} of the graded maximal ideal. This terminology conforms with that of \cite{gasharov1998green} and \cite{gasharov2002resolutions}.

\begin{definition}\label{def: restrictedPower} Let $\m = (x_1,\dots, x_n)$ be the graded maximal ideal in a polynomial ring $S$ over a field $k$. For any vector $\uu\in \NN^n$, define the \emph{restricted power} of $\m$ to be $\m^\dd(\leq \uu)$ to be the ideal generated by $(\xx^\aa \mid \sum_i^n a_i = d \text{ and } a_i\leq u_i \text{ for all } i\in [n])$. 
\end{definition}

\begin{setup}\label{set: generalizedPols}
Let $S = k[x_1,\dots, x_n]$ be a polynomial over a field $k$. Denote by $\mathfrak{m} = (x_1,\dots, x_n)$ the graded maximal ideal of $S$. For $\uu\in \NN^n$, let $\m^d(\leq \uu)$ be the restricted power of the graded maximal ideal as in Definition \ref{def: restrictedPower}.
Let $\Xv_i = \{x_{i1},\dots, x_{i d}\}$ be a set of variables, and let $\widetilde S = k[\Xv_1,\dots, \Xv_n]$ be a polynomial ring in the union of all these variables. Let $\cB_d(\leq u_i)$ be the \emph{truncated Boolean poset} on $[d]$ with elements of rank at most $u_i$ in $\cB_d$.

Let $\cH_n^d$ be the hypersimplicial complex from Definition \ref{def: hypersimplex}, and let $\cH_n^d(\leq \uu)$ be the induced subcomplex with cells $\{ C_{\aa,J} \mid \mdeg(C_{\aa,J}) \leq \uu\}$.
Denote by $\sk_0(\cH^d_n (\leq \uu))$ and $\sk_1(\cH^d_n (\leq \uu))$ the $0$-skeleton and $1$-skeleton of $\cH_n^d(\leq \uu)$, respectively. We also use the notation $\Delta^\ZZ_{\leq \uu}(n,d)$ for $\sk_0(\cH^d_n (\leq \uu))$.

Let $L^a_{b,\leq \uu}(F)$ be the Schur module defined in Setup \ref{set:Lcomplexsetup} restricted to multidegrees $\leq \uu$, where the multidegree of an element $f_J \otimes f^\alpha \in \bigwedge F^a \otimes S_b(F)$ is defined to be $\alpha + \sum_{j\in J} \epsilon_j$, where $\epsilon_j$ is the $j$'th unit vector in $\NN^n$.
\end{setup}

\begin{prop}\label{prop: boundedPowersCellRes} Adopt notation and hypotheses of Setup \ref{set: generalizedPols}. Then:
\begin{enumerate}
    \item The induced subcomplex $\cH^d_n(\leq \uu)$ supports a polyhedral cellular resolution of $\m^d(\leq \uu)$.
    \item The induced subcomplex $\widetilde{\cH}^d_n(\leq\uu)$ supports a minimal CW cellular resolution of $\m^d(\leq \uu)$ which is isomorphic to a subcomplex of the $L$-complex.
\end{enumerate}
In particular, $\m^d(\leq \uu)$ has a linear minimal free resolution.
\end{prop}

\begin{proof}
Apply Proposition \ref{prop: subresolution}.
\end{proof}

\begin{cor}\label{cor: squarefreeCellRes} Adopt notation and hypotheses of Setup \ref{set: boundedGens} and let $\one = (1,1,\dots, 1)\in \NN^n$. The ideal generated by all squarefree monomials of a given degree in $S$ has a non-minimal resolution supported on the polyhedral cell complex $\cH_n^d(\leq \one)$, and it has a minimal free resolution supported on the CW-complex $\widetilde{\cH}_n^d(\leq \one)$, the restriction of the Morse complex from Proposition \ref{prop: cellularLcplx}.
\end{cor}

Moreover, one can extend the characterizations of polarizations of powers of the graded maximal ideal in Theorem \ref{thm: spanningTree} to restricted powers of the maximal ideal. Observe that all the definitions in Section \ref{sec: pols} work in this context, exchanging $\Delta^\ZZ(n,d)$ with $\sk_0(\cH_n^d(\leq \uu))$ and exchanging $\cB_d$ with $\cB_d(\leq u_i)$ as required.

\begin{theorem}\label{thm: generalizedPols} Adopt notation and hypotheses of Setup \ref{set: boundedGens}. Let $\chi = \{X_i\}_{i\in [n]}$ denote a set of rank-preserving isotone maps
$$
X_i: (\Delta^{\ZZ}_{\leq \uu}(n,d), \leq_i) \ra \cB_d(\leq u_i)
$$
as in Construction \ref{const: pols}.
Denote by $\tabs(\chi)$ be the set of tableaux in $\bigwedge^2 \FF \otimes S_{d-1}(F)$ associated to the linear syzygy edges in $\sk_1(\cH_n^d (\leq \uu))$ after applying the isotone maps in $\chi$ to its vertices.
Then the following are equivalent:
\begin{enumerate}
    \item The elements of $\tabs(\chi)$ span the module $L_{d, \leq \uu}^1(F)$.
    \item For every $\cc\in\Delta^\ZZ_{\leq \uu}(n,d+1)$, $\LS(\cc)$ contains a spanning tree of the complete down-graph $D(\cc)$.
    \item The set of isotone maps $X_1,\dots, X_n$ determine a polarization of $\m^d(\leq \uu)$.
\end{enumerate}
\end{theorem}

\begin{proof}
\noindent (1) $\implies$ (2): The proof is identical to that of Theorem \ref{thm: spanningTree}.

\noindent (2) $\implies$ (3): This follows from Proposition \ref{prop: restrictedPols}.

\noindent (3) $\implies$ (1):
Suppose the set of isotone maps $\chi = \{X_i\}_{i\in [n]}$ determine a polarization $\widetilde{\m^d(\leq \uu)}$ of $\m^d(\leq\uu)$. Let $\widetilde \FF$ be a (not necessarily minimal) free resolution of $\widetilde{\m^d(\leq \uu)}$ with linear syzygies corresponding to the set of linear syzygy edges induced by $\chi$. Let $\FF = \widetilde \FF \otimes \widetilde S / \sigma$, be the depolarization of $\widetilde \FF$, where $\sigma$ is a regular sequence of variable differences. Then $\FF$ is a free resolution of $\m^d(\leq\uu)$ with linear syzygies which are in bijection with $\tabs(\chi)$. Since $\FF$ must be a free resolution of $\m^d(\leq\uu)$, $\tabs(\chi)$ must span $L^1_{d,\leq\uu}(F)$.

\end{proof}

\section*{Acknowledgments} We would like to thank Gunnar Fl\o ystad, Benjamin Smith, and the anonymous referees for helpful feedback on earlier drafts of this paper. The first author was partially supported by the NSF GRFP under Grant No. DGE-1650441. 
\bibliographystyle{amsplain}
\bibliography{biblio}

\providecommand{\bysame}{\leavevmode\hbox to3em{\hrulefill}\thinspace}
\providecommand{\MR}{\relax\ifhmode\unskip\space\fi MR }
\providecommand{\MRhref}[2]{%
  \href{http://www.ams.org/mathscinet-getitem?mr=#1}{#2}
}
\providecommand{\href}[2]{#2}
\begin{thebibliography}{10}

\bibitem{polarizations}
Ayah Almousa, Gunnar Fl{\o}ystad, and Henning Lohne, \emph{Polarizations of
  powers of graded maximal ideals}, Journal of Pure and Applied Algebra
  \textbf{226} (2022), no.~5, 106924.

\bibitem{batzies2002discrete}
Ekkehard Batzies, \emph{Discrete {M}orse theory for cellular resolutions}, J.
  reine angew. Math \textbf{543} (2002), 147--168.

\bibitem{bayer1998monomial}
Dave Bayer, Irena Peeva, and Bernd Sturmfels, \emph{Monomial resolutions},
  Mathematical Research Letters \textbf{5} (1998), no.~1, 31--46.

\bibitem{buchsbaum1975generic}
David~A Buchsbaum and David Eisenbud, \emph{Generic free resolutions and a
  family of generically perfect ideals}, Advances in Mathematics \textbf{18}
  (1975), no.~3, 245--301.

\bibitem{clark2012minimal}
Timothy~BP Clark, \emph{A minimal poset resolution of stable ideals}, Progress
  in commutative algebra 1, 2012, pp.~143--166.

\bibitem{el2014artinian}
Sabine El~Khoury and Andrew~R Kustin, \emph{Artinian {G}orenstein algebras with
  linear resolutions}, Journal of Algebra \textbf{420} (2014), 402--474.

\bibitem{forman1998morse}
Robin Forman, \emph{Morse theory for cell complexes}, Advances in Mathematics
  \textbf{134} (1998), 90--145.

\bibitem{forman2002user}
\bysame, \emph{A user's guide to discrete {M}orse theory.}, S{\'e}minaire
  Lotharingien de Combinatoire [electronic only] \textbf{48} (2002), B48c--35.

\bibitem{gasharov1998green}
Vesselin Gasharov, \emph{Green and {G}otzmann theorems for polynomial rings
  with restricted powers of the variables}, Journal of Pure and Applied Algebra
  \textbf{130} (1998), no.~2, 113--118.

\bibitem{gasharov2002resolutions}
Vesselin Gasharov, Takayuki Hibi, and Irena Peeva, \emph{Resolutions of
  a-stable ideals}, Journal of Algebra \textbf{254} (2002), no.~2, 375--394.

\bibitem{hartshorne1966connectedness}
Robin Hartshorne, \emph{Connectedness of the {H}ilbert scheme}, Publications
  Math{\'e}matiques de l'Institut des Hautes {\'E}tudes Scientifiques
  \textbf{29} (1966), no.~1, 7--48.

\bibitem{mermin2010eliaiiou}
Jeffrey Mermin, \emph{The {E}liahou-{K}ervaire resolution is cellular}, Journal
  of Commutative Algebra \textbf{2} (2010), no.~1, 55--78.

\bibitem{peeva2010graded}
Irena Peeva, \emph{Graded syzygies}, vol.~14, Springer Science \& Business
  Media, 2010.

\bibitem{peeva2011frames}
Irena Peeva and Mauricio Velasco, \emph{Frames and degenerations of monomial
  resolutions}, Transactions of the American Mathematical Society (2011),
  2029--2046.

\bibitem{taylor66}
Diana~Kahn Taylor, \emph{Ideals generated by monomials in an r-sequence}, Ph.D.
  thesis, University of Chicago, Department of Mathematics, 1966.

\bibitem{velasco2008minimal}
Mauricio Velasco, \emph{Minimal free resolutions that are not supported by a
  cw-complex}, Journal of Algebra \textbf{319} (2008), no.~1, 102--114.

\bibitem{weyman2003}
Jerzy Weyman, \emph{Cohomology of vector bundles and syzygies}, vol. 149,
  Cambridge University Press, 2003.

\end{thebibliography}
\addcontentsline{toc}{section}{Bibliography}

\end{document}